\documentclass[preprint,12pt]{elsarticle}

\usepackage{multirow}
\usepackage{float}
\usepackage{algorithm}
\usepackage{algorithmic}
\usepackage{graphicx}
\usepackage{amssymb,amsmath}
\usepackage{subfigure}

\usepackage{float}
\biboptions{sort&compress}

\newtheorem{example}{Example}[subsection]

\newproof{proof}{Proof}
\numberwithin{equation}{section}

\usepackage{lineno}
\usepackage{xcolor}

\newcommand{\vv}{\overrightarrow{V}}

\newcommand{\uAve}{\bar{u}}

\newcommand{\vW}{\textbf{W}}
\journal{Journal name}
\UseRawInputEncoding
\begin{document}

\begin{frontmatter}



\title{Cell-average based neural network method for high dimensional parabolic differential equations}

 \author[label1]{Hong Zhang}
 \ead{hzcomath@163.com}
 \author[label1,label2]{Hongying Huang\corref{cor1}}
 \ead{huanghy@lsec.cc.ac.cn}
 \author[label3]{Jue Yan}
 \ead{jyan@iastate.edu}
 \address[label1]{Zhejiang Ocean University,
          Zhoushan, Zhejiang 316000, P.R. China.}
 \address[label2]{Key Laboratory of Oceanographic Big Data Mining \& Application of Zhejiang Province, Zhoushan, Zhejiang 316000, P.R. China}
 \address[label3]{Department of Mathematics, Iowa State University, Ames, IA 50011, USA}
%\author[label1]{Hong Zhang}
%%\ead{cxqiu@iastate.edu}

\fntext[label2]{The research of Huang was supported by National Natural Science Foundation of China under grant NO.11771398.  }
\fntext[label3]{Research work of Yan is supported by National Science Foundation grant DMS-1620335 and Simons Foundation grant 637716.}
\cortext[cor1]{corresponding author}

\begin{abstract}
In this paper, we introduce cell-average based neural network (CANN) method to solve high-dimensional parabolic partial differential equations. The method is based on the integral or weak formulation of partial differential equations.
A feedforward network is considered to train the solution average of cells in neighboring time. Initial values and approximate solution at $t=\Delta t$ obtained by high order numerical method are taken as the inputs and outputs of network, respectively. We use supervised training combined with a simple backpropagation algorithm to train the network parameters.
We  find that the neural network  has been trained to optimality for high-dimensional problems, the CFL condition is  not strictly limited for CANN method and the trained network is used to solve the same problem with different initial values. For the high-dimensional parabolic equations, the convergence is observed and the errors are shown related to spatial mesh size but independent of time step size.
\end{abstract}

\end{frontmatter}

\section{Introduction}
\label{sec1}

In this paper, we consider the following high-dimensional parabolic partial differential equations:
\begin{equation}\label{eq:model}
    u_t+\beta \nabla \cdot f(u)=\mu \Delta u+g(u),~~~(\mathbf{x},t)\in \Omega \times R^{+},
\end{equation}
where  $t$ and $\mathbf{x}$ denote the time and spatial variables respectively, and $\Omega \in R^d$ is the spatial domain and $d\geq 2$.  We develop cell-average based neural network (CANN) method \cite{qiu2021cell} to solve high-dimensional parabolic equation (\ref{eq:model}). The
method is based on the integral or weak formulation of partial differential equations similar to finite volume method.

Machine learning with neural networks have achieved great success in image classification, pattern recognition and cognitive science \cite{bengio2009learning,beck2019machine,goodfellow2016deep,lecun2015deep} for the last two decades. In the last few years, machine learning based on neural networks have also been explored to solve partial differential equations \cite{rackauckas2020universal,chen2018neural,long2019pde,ruthotto2020deep,he2019mgnet}.
These methods are classified into two categories:  a better solver may be obtained or it can assist to enhance the performance of current numerical methods

One class is to combine network with classical numerical methods for performance improvement.
The work of \cite{ray2018artificial} applied neural networks as trouble-cell indicator. Deep reinforcement network in \cite{wang2019learning} is explored to estimate the weights and enhance WENO schemes performance. Neural network \cite{discacciati2020controlling} is applied for identifying suitable amount of artificial viscosity added.

The other  is to employ  neural network to directly approximate the solution of PDEs \cite{lagaris1998artificial,rudd2015constrained,sirignano2018dgm,raissi2017physics,raissi2019physics,dwivedi2020physics,lu2021deepxde,jin2021nsfnets,shin2020convergence,laakmann2021efficient,cai2021least,cai2021least1,qiu2021cell}. With $x$ and $t$ as network input vector, such methods have the advantage of automatic differentiation, mesh free and can be applied to solve many types of equations.% Approximation properties of a class of artificial neural networks are established by \cite{cybenko1989approximation,hornik1990universal,barron1993universal,leshno1993multilayer,pinkus1999approximation}.
The popular PINN methods in \cite{raissi2017physics,raissi2019physics} is
proposed. An extreme and distributed network to improve the PINN efficiency on larger domain in \cite{dwivedi2020physics} is considered.  Application to incompressible Navier-Stokes equations is studied in \cite{jin2021nsfnets}, which could overcome some of the aforementioned limitations for simulating incompressible laminar and turbulent flows. We refer to \cite{shin2020convergence,laakmann2021efficient} for the development of PINN convergence.
In \cite{Weinan2018Ritz}, a deep Ritz method  is proposed to solve the class of PDEs that can be reformulated as equivalent energy minimization problems. The constraint due to boundary condition is added to the energy as a penalty term in \cite{Weinan2018Ritz}. The weak adversarial network that could convert the problem of finding the weak solution of PDEs into an operator norm minimization problem induced from the weak formulation, is proposed by \cite{zang2020weak}. In \cite{cai2021least} a discretization of an equivalent least-squares formulation in the set of neural network functions with the ReLU activation function is used to solve the advection-reaction equations. Very recently, the cell-average based neural network method which is closely related to finite volume scheme has been proposed in \cite{qiu2021cell}. This method can relief from  CFL restrictions for  explicit difference scheme and can adapt to any time size for neural networks solver. Motivated by simple of the method and its successful numerical experiments, we further study CANN method for more higher dimensional parabolic problems.

High-dimensional partial differential equations are used in physics, engineering, and finance. Their numerical solution has been a longstanding challenge.
Classical numerical methods become infeasible in high-dimensions due to the explosion in the number of grid points and the demand for reduced time step size. If there are $d$ space dimensions and 1 time dimension, the number of the mesh elements is of  $\mathcal{O}(1/h^{d+1})$.  This quickly becomes computationally intractable when the dimension $d$ becomes moderately large. The classical approach, by building functions using polynomials, wavelets or other basis functions, is bound to run into the curse of dimensionality problem. But the neural network model has shown remarkable success in  solving high-dimensional deterministic differential equations \cite{beck2019machine,chan2019machine,pham2021neural,hutzenthaler2020proof,lye2020deep,wu2020data,qin2021data}.

In this paper, we extend  CANN method to solve high-dimensional parabolic problem. The rest of article is organized follows. In Section \ref{S:2.1}, We introduce CANN method and highlight its connection to finite volume scheme. In Section \ref{S:2.2} and Section \ref{S:2.3} , we present the network structure and loss function in training process. Some numerical examples are shown in Section \ref{sec:numerical examples}. Finally some conclusions are given in Section \ref{S:4}.

\section{Cell-average based neural network method}
\label{sec:neural network}

\subsection{Problem setup, motivation}
\label{S:2.1}

In this section we develop  cell-average based neural network  method solving high-dimensional partial differential equations
\begin{equation}\label{eq:PDE}
\begin{aligned}
     u_t &= \mathcal{L}(u),~~~(\mathbf{x},t) \in \Omega \times R^{+},\\
     u &= u|_{t=0},\\
     u &= u_0, ~  \text{on}~ {\partial \Omega},
\end{aligned}
\end{equation}
where
differential operator $\mathcal{L}$ is introduced to represent a generic high-dimensional differential operator. This setup encapsulates a wide range of problems in mathematical physics including conservation laws, diffusion processes and advection-diffusion-reaction system. For example, we have $\mathcal{L}(u)=u_{xx}+u_{yy}$ for two dimensional heat equation.

\begin{figure}
    \centering
    \includegraphics[width=5.5cm]{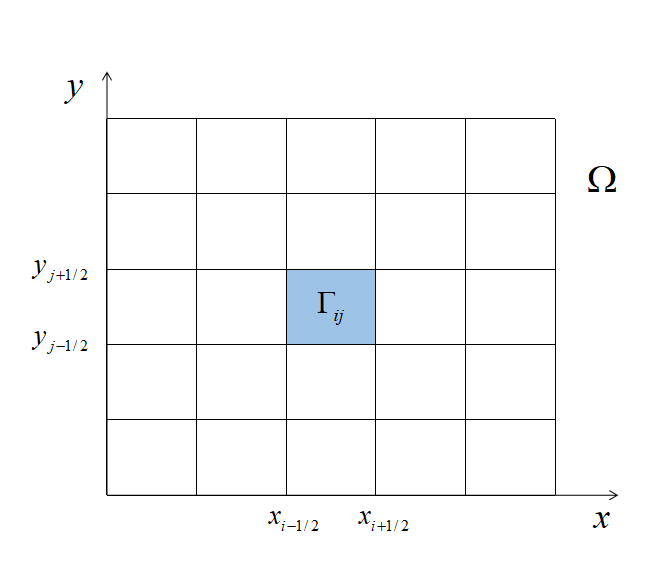}
    \caption{Partition of cell-average based neural network method}
    \label{partition}
\end{figure}

For simplicity of presentation, we focus on two dimensional problems and $\Omega=[a,b]\times [c,d]$. The CANN method is motivated by finite volume method.  The essential step is to divide $\Omega$ into many control volumes and approximate the integral conservation law on each of the control volumes. Divide the domain $\Omega$ into $I \times J$ cells and let $\Delta x=(b-a)/I, \Delta y=(d-c)/J$ and  $x_{1/2}=a,\cdots, x_{i+1/2}=a+(i-1) \Delta x,\cdots, x_{I+1/2}=b$ and $ y_{1/2}=c, \cdots,  y_{j+1/2}=c+(j-1)\Delta y,\cdots, y_{J+1/2}=d$. Let $\Gamma_{ij}=[x_{i-1/2} , x_{i+1/2}] \times [y_{j-1/2}, y_{j+1/2}]=I_i \times J_j$ denote the cell in column $i$ and row $j$ as shown in Fig. \ref{partition} and  the area of the cell $\Gamma_{ij}$ is $\Delta S=\Delta x \times \Delta y$. Let $\Delta t$ be time size and $t_n=n\times \Delta t, n=0,1,\cdots$.

Integrating  partial differential Eq. (\ref{eq:PDE}) over the  computational cell $\Gamma_{ij}$ and time interval $[t_n, t_{n+1}]$, we have
\begin{equation}\label{eq:PDE-integration}
\int^{t_{n+1}}_{t_n}\int_{\Gamma_{ij}}  u_t~dxdydt=\int^{t_{n+1}}_{t_n}\int_{\Gamma_{ij}} \mathcal{L}(u)~ dxdydt.
\end{equation}
 Let $\bar{u}_{ij}(t)$ denote the cell average of $u$  in $\Gamma_{ij}$, that is
$\bar{u}_{ij}(t)=\frac{1}{\Delta S}\int_{\Gamma_{ij}} u(x,y,t) ~dxdy$, then
Eq. (\ref{eq:PDE-integration}) can be integrated out as
\begin{equation}\label{eq:PDE-integral-format}
\bar{u}_{ij}(t_{n+1}) -\bar{u}_{ij}(t_{n})=\frac{1}{\Delta S}\int^{t_{n+1}}_{t_n}\int_{\Gamma_{ij}} ~ \mathcal{L}(u)~ dxdydt.
\end{equation}
Referring to \cite{qiu2021cell}, we use a simple fully connected network $\mathcal{N}(\cdot; \Theta)$ to estimate the terms at the right of (\ref{eq:PDE-integral-format})
\begin{equation}\label{nn:network-goal}
    \mathcal{N}(\cdot; \Theta)\approx \frac{1}{\Delta S}\int^{t_{n+1}}_{t_n}\int_{\Gamma_{ij}} ~ \mathcal{L}(u)~ dxdydt,
\end{equation}
where $\Theta$ denotes the network parameter set of all weight matrices and biases.

According to (\ref{eq:PDE-integral-format}) and (\ref{nn:network-goal}), we could use the solution average $\left\{\bar{u}^n_{ij}\right\}$ at time level $t_n$ to approximate the solution average $\bar{u}^{n+1}_{ij}$ at the next time level $t_{n+1}$ by the following neural network
\begin{equation}\label{nn:neural-network}
    \bar{v}_{ij}^{out} = \bar{v}_{ij}^{in}+\mathcal{N}(\vv_{ij}^{in}; \Theta).
\end{equation}

Considering $\bar{v}_{ij}^{in}=\bar{u}^{n}_{ij}$ and comparing (\ref{nn:neural-network}) with the integral form (\ref{eq:PDE-integral-format}) of the PDEs, we have
$$
\bar{v}_{ij}^{out}\approx \bar{u}^{n+1}_{ij}.
$$
\begin{figure}
    \centering
    \includegraphics[width=9.5cm]{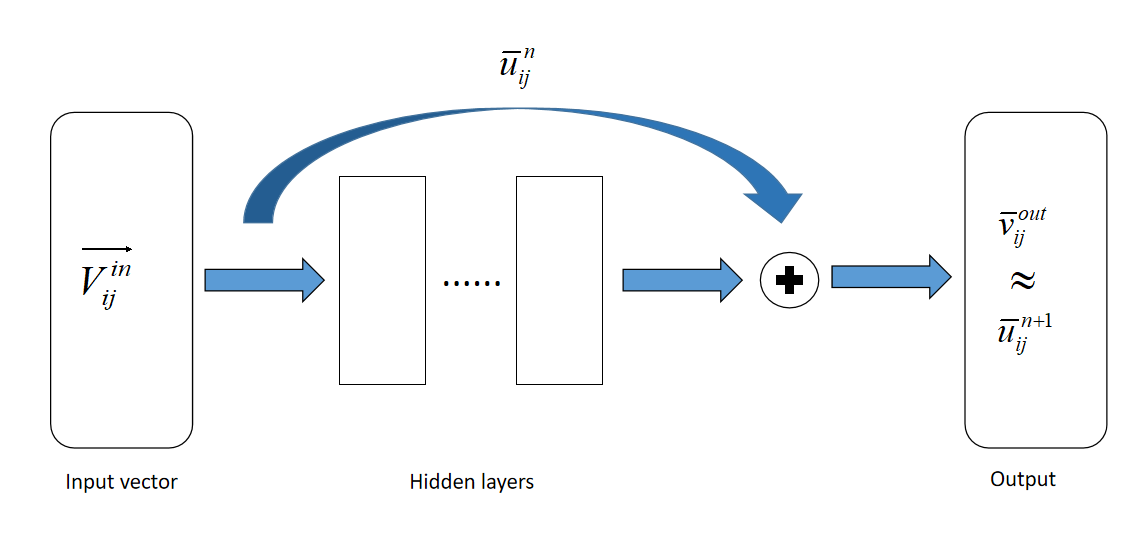}
    \caption{Illustration of cell-average based neural network method}
    \label{illustration}
\end{figure}

\subsection{Network architecture}
\label{S:2.2}

\begin{figure}
    \centering
    \includegraphics[width=5.5cm]{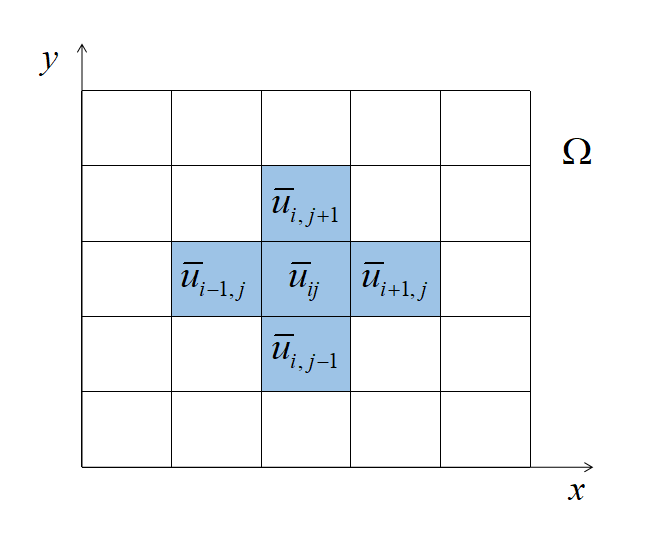}
    \caption{Input vectors with common edges}
    \label{five points}
\end{figure}

As shown in Fig. \ref{illustration}, the vector $\vv_{ij}^{in}$ as the input vector is the base of the network structure.
We consider two types of  input vectors of neural network. One is similar to five-point difference scheme of finite difference method
\begin{equation}\label{nn:input-vector-stencil}
    \vv_{ij}^{in}=\Big[\bar{u}^n_{i-1,j},  \bar{u}^n_{i+1,j}, \bar{u}^n_{ij}, \bar{u}^n_{i,j+1},\bar{u}^n_{i,j-1} \Big]^T,
\end{equation}
where the cells have the common edges with the cell $\Gamma_{ij}$ as shown in Fig. \ref{five points}.

\begin{figure}
    \centering
    \includegraphics[width=5.5cm]{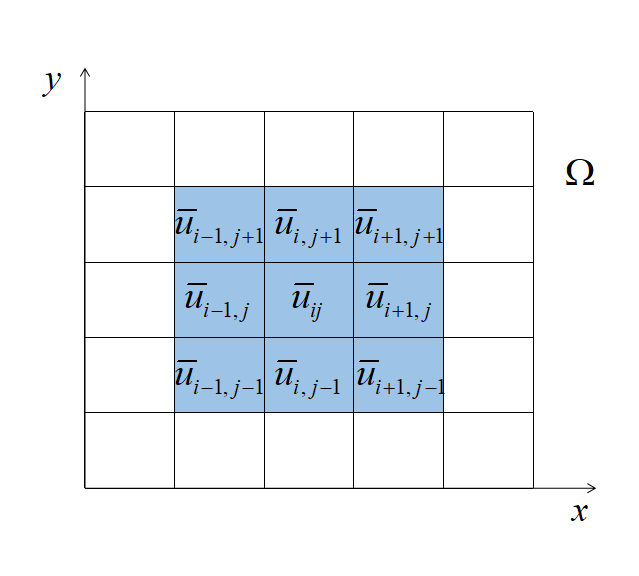}
    \caption{Input vectors with common vertices}
    \label{nine points}
\end{figure}

The other is similar to nine-point difference scheme, which share the same vertices with the cell $\Gamma_{ij}$ in Fig. \ref{nine points}
\begin{equation}\label{nn:co-vertex-input-stencil}
    \vv_{ij}^{in} = \Big[ \uAve^n_{i-1,j+1},\uAve^n_{i,j+1},\uAve^n_{i+1,j+1}, \uAve^n_{i-1,j},\uAve^n_{ij},\uAve^n_{i+1,j},\uAve^n_{i-1,j-1},\uAve^n_{i,j-1},\uAve^n_{i+1,j-1}\Big]^T.
\end{equation}

The scheme (\ref{nn:input-vector-stencil}) or (\ref{nn:co-vertex-input-stencil})  may influence the effectiveness  of the neural network method approximating the solution average $\bar{u}_{ij}^{n+1}$ at the next time level.

%%介绍神经网络的结构
We use a simple feedforward neural network with $M$ layers and $M\geq 3$. Set the input and output as  the first and last layer, respectively, then the hidden layer is $M-2$ layers. In Fig. \ref{illustration}, we show the network structure of  CANN method.  Then taking $n_i ~(i=1,\cdots, M)$ as the number of neurons in $i$-th layer. The number of neurons in the first layer depends on the number of cells in the input vector (\ref{nn:input-vector-stencil}) or (\ref{nn:co-vertex-input-stencil}), and is defined as $n_1$. Meanwhile, we can find $n_{M}=1$ from (\ref{nn:neural-network}). This neural network can be considered as a function $\mathcal{N}: R^{n_1} \to R^{n_M}$.

Learning data are collected in the form of pairs. Each pair refers to the solution averages at
 two neighboring time level, denoted as
\begin{equation}\label{nn:learning-data-set}
S =\left\{\left(\vv_{ij}^{in},\bar{u}_{ij}^{n+1}\right),\,\,i=1, \cdots, I , j=1, \cdots, J\right\}_{n=0}^p,
\end{equation}
where $\vv_{ij}^{in}$ is five-point scheme (\ref{nn:input-vector-stencil}) or nine-point scheme (\ref{nn:co-vertex-input-stencil}). The solution averages $\bar{u}_{ij}^n$ and $\bar{u}_{ij}^{n+1}$ are obtained from high order numerical method for the PDEs, such as  discontinuous Galerkin method or finite difference method or measured from the lab for application problems. When exact solution for partial differential equations, one time level solution averages in the learning data $S$ corresponding to $(t_0,t_1)$ or $p=0$ in (\ref{nn:learning-data-set}) is sufficient for obtaining an effective neural network(see the results in Example \ref{ex1}-\ref{ex8}).

And so much experience has suggested the smoothness of the activation function $\sigma_{i}$ , where $i$ is the layer of all network, plays a important role in the accuracy of the algorithm. To balance simplicity and accuracy, we have decided to use $\tanh(x)$. Specially, $\sigma_i=\tanh(x)$ is applied between all layers, except to the output layer for which we have $\sigma_M=x$. The full $M$-layer network can now be expressed as:
\begin{equation}
        \bar{v}_{ij}^{out} = \bar{v}_{ij}^{in}+\mathcal{N}(\vv_{ij}^{in}; \Theta)
=\bar{v}_{ij}^{in}+(\sigma_M\circ \vW_{M-1})\circ \cdots \circ (\sigma_2\circ \vW_1)(\vv_{ij}^{in}), \label{neuralNetwork-structure}
\end{equation}
where $\circ$ stands for operator composition. We have $\vW_i$ denoting the linear transformation operator or the weight matrix connecting the neurons from $i$-th layer to $(i+1)$-th layer. Furthermore, we proposed that  CANN method, once well trained, would be implemented as a regular explicit finite volume scheme.

\subsection{Training data and testing data}
\label{S:2.3}

We present three options for training set and the corresponding testing set.

\begin{enumerate}\item For given initial-boundary value problem of partial differential equation, we take
the solution average data from $t_0$ to $t_1$ as the training set and then use the well trained network solver to evolve numerical
solution forward in time. Accumulation errors in time can be small and the network solver is still stable.

\item We randomly select $75\%$ of solution average data (\ref{nn:learning-data-set}) with $p=0$ from $t_0$ to $t_1$ as the training set and the other $25\%$ is used as the test data. Then we obtain the approximate solution at $t_n, n\geq 2$ time level from the trained network.

\item For given initial-boundary value problem of partial differential equation, we take
the solution average data from $t_0$ to $t_1$ as training set and  then use the well trained network solver to solve the same PDEs with
different initial data. For example, we consider the heat equation with periodic boundary values and initial value $u(x,y,0)=\sin(x+y)$. Firstly, we use numerical method for this problem to obtain the cell average $\bar{u}_{ij}^1$ at $t_1$ time level. Secondly, we apply training set with (\ref{nn:learning-data-set}) to train the network. Finally, we use the well trained network to solve the same PDE with initial value $u(x,y,0)=\cos(x+y)$ or $u(x,y,0)=\cos(x+y+\pi/3)$.
\end{enumerate}

\subsection{Training process}
\label{S:2.4}

In this section, we train the neural network to get the optimal parameter $\Theta^{*}$ in (\ref{neuralNetwork-structure}).
Our goal is to use the trained parameter to accurately approximate the solution average evolution $\uAve^n_{ij}\to \uAve^{n+1}_{ij}$. First, we apply (\ref{neuralNetwork-structure}) with  $\bar{v}^{in}_{ij}=\uAve_{ij}^n$ to obtain the network output $\bar{v}_{ij}^{out}$. Then it is compared with the target output $\uAve_{ij}^{n+1}$. The square loss function donated as
\begin{equation}\label{nn:L2-error-square}
   L_{t_{n+1}}(\Theta)=\sum^I_{i=1}\sum^J_{j=1}(\bar{v}_{ij}^{out}-\uAve_{ij}^{n+1})^2\Delta S,
\end{equation}
for all $i=1,\cdots,I$ and $j=1,\cdots,J$.  Notice that in the input vector (\ref{nn:input-vector-stencil}) or (\ref{nn:co-vertex-input-stencil}), the solution average of the boundary cell is closely related to the squared loss function.  For CANN method (\ref{nn:neural-network}), we assign ghost cell values and apply boundary conditions as a regular finite volume method. In this paper, to be more efficient, we only consider Dirichlet or periodic boundary conditions. Finally, we loop among the data set $S$ in (\ref{nn:learning-data-set}) to minimize the error by (\ref{nn:L2-error-square}). The stochastic or approximated gradient descent method is applied for a single data pair. Given tolerance $\epsilon$, we minimize the error (\ref{nn:L2-error-square})  to acquire the optimal parameter set $\Theta^{*}$ .

\begin{algorithm}[htb]
\caption{Algorithm for CANN method}
\label{alg:cann_alg}
\begin{algorithmic}[1]
\STATE Generate $\uAve_{ij}^{n}$ and $\uAve_{ij}^{n+1}$ by high order numerical method where $i=1,2,\ldots,I,~ j=1,2,\ldots, J$
%\STATE Use a simple fully connected neural network $\mathcal{N}(\cdot; \Theta)$
\STATE Generate the input vector $\vv_{ij}^{in}$ by (\ref{nn:input-vector-stencil}) or (\ref{nn:co-vertex-input-stencil}), and target data $\uAve_{ij}^{n+1}$, where $i=1,2,\ldots,I, j=1,2,\ldots, J$
\STATE  Use a simple fully connected neural network $\mathcal{N}(\cdot; \Theta)$, and $\bar{v}_{ij}^{out} = \uAve_{ij}^{in}+\mathcal{N}(\vv_{ij}^{in}; \Theta)$
\STATE Compute the minimizer of the expected square loss function by stochastic gradient descent method
\begin{equation}\label{alg:L2_loss}
\left\{
\begin{array}{l}
L_{t_{n+1}}(\Theta)= \sum^I_{i=1}\sum^J_{j=1}(\bar{v}_{ij}^{out}-\uAve_{ij}^{n+1})^2\Delta S\\
\Theta^{*} \in \arg ~ \min~ L_{t_{n+1}}(\Theta).
\end{array}
\right.
\end{equation}
\STATE Update: $\bar{v}_{ij}^{out}.$

\end{algorithmic}
\end{algorithm}

\section{Numerical examples} \label{sec:numerical examples}

In this section, we apply CANN method based on the loss function defined in (\ref{nn:L2-error-square}) to solve a series of high-dimensional parabolic problems. Firstly, we use discontinuous Galerkin method for two-dimensional problem or finite difference method for higher-dimensional problems to obtain the approximate solution with high order accuracy at $t=t_1=\Delta t$. Secondly, we take the initial value as input vectors and the approximate solution at $t=\Delta t$ as output vectors and then train neural network. Finally, we compute the approximate solution at $t=t_{n+1}$  for given solution at $t=t_n$ with the trained optimal network. In order to check out the accuracy and capability of CANN method, we vary the number of layer and neurons per layer.

We compute the error and order at the final time which define as $T$. For simplicity of presentation, we always choose $\Delta x=\Delta y$ with taking the same space step size. We consider two error norms:

\begin{equation}\label{nn:L2-norm-error}
  e_{L_2}(T)=\sqrt{\sum^I_{i=1}\sum^J_{j=1}(\bar{v}_{ij}(T)-\uAve_{ij}(T))^2\Delta S}
\end{equation}
\begin{equation}\label{nn:Linfty-norm-error}
   e_{L_{\infty}}(T)=\max_{i}\max_{j}|\bar{v}_{ij}(T)-\uAve_{ij}(T)| \end{equation}

After we train the networks, we may give errors and orders in (\ref{nn:L2-norm-error}) and (\ref{nn:Linfty-norm-error}) with the optimal parameter $\Theta^*$. Again we emphasis that we should follow the principle listed in Section \ref{S:2.2} to choose suitable network input vector. We also mention that squared $L_2$ errors of (\ref{nn:L2-error-square}) is used as training break condition, with which well trained network errors are around $10^{-6}$ or smaller.

 We start with linear diffusion equations and then we move on more complicated nonlinear diffusion problem.

%%%%%%%%%%%%%%%%%%%%%%%%%%%%%%%%%%%%%%%%%%%%%%%%%%%%%%%%%%%%%%%%%%
\subsection{2D diffusion equations}
\label{S:3.1}

\begin{example}\label{ex1}{\bf\emph {Heat equation}} \end{example} Let $\Omega =(0,2\pi)\times (0,2\pi)$, consider the heat equation: \begin{equation}\label{eq:3-1}
u_t  = u_{xx}+u_{yy}, \ \text{in} \ \Omega \times (0,T), \quad \text{with}\ \
u(x,y,0) =\sin(x+y)，
\end{equation}
with periodic boundary condition.  We use discontinuous Galerkin method for this problem to obtain the approximate solution average at the time  $t_1$  as the training target and then get well trained network.  Finally, we generalized  the trained network to solve the same PDE with the different initial value $u(x,y,0)=\cos(x+y)$ and $u(x,y,0)=\cos(x+y+\pi/3)$. Let the number of iteration $K=10^{5}$ and use 1 hidden layer and 10 neurons per layer. We use final time level $T=\pi$ to compute all errors.

\vspace{.1in}
\noindent {\bf Case I: test accuracy of five-point and nine-point scheme. }

\begin{figure}[htbp]
   \centering \subfigure[Example \ref{ex1}: Heat equation]{
   \includegraphics[width=6.5cm]{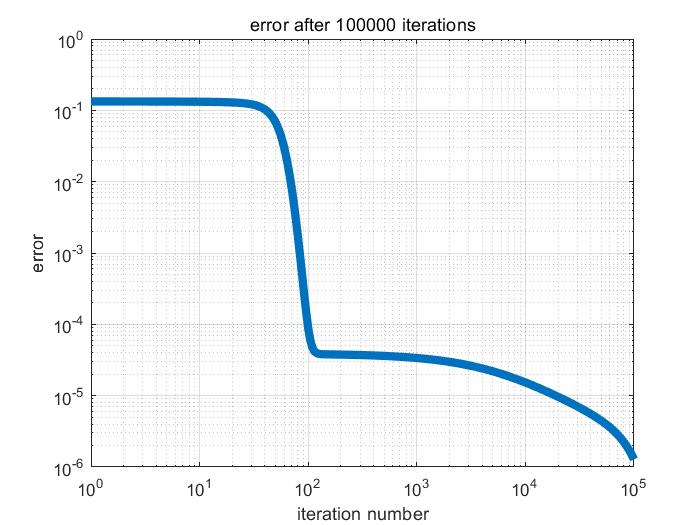} } \quad \subfigure[Example
   \ref{ex2}: Convection diffusion equation]{
   \includegraphics[width=6.5cm]{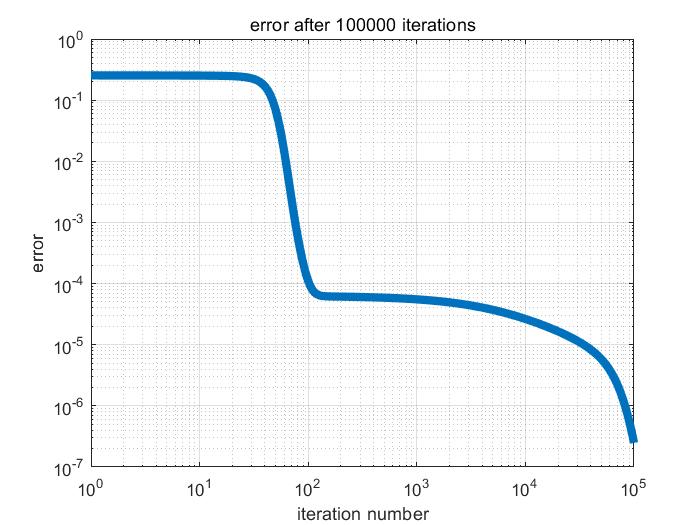} } \quad \subfigure[Example
   \ref{ex3}: Anisotropic diffusion equation]{
   \includegraphics[width=6.5cm]{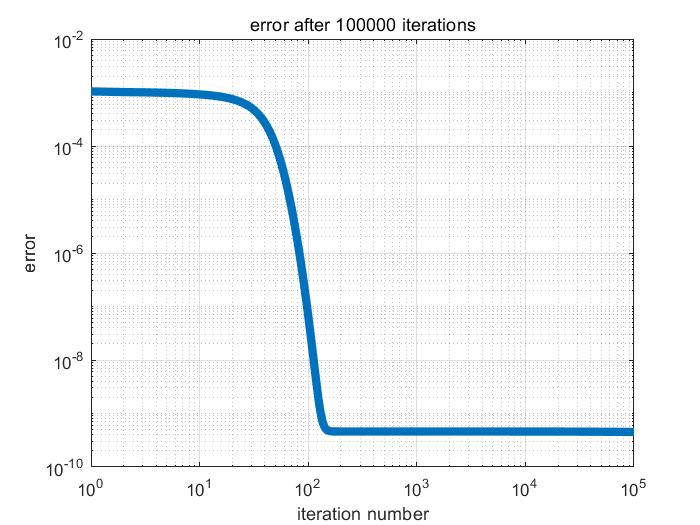} } \quad \subfigure[Example
   \ref{ex4}: Navier-Stokes equation]{
   \includegraphics[width=6.5cm]{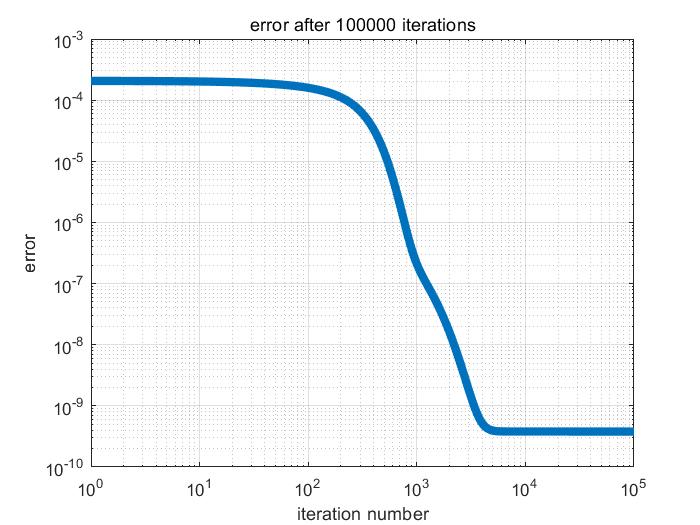} } \caption{Neural network squared $L_2$ training errors
   (\ref{nn:L2-error-square})}\label{1-4errors}
   \end{figure}

For first case , we test the accuracy of a standard highly dimensional CANN method with five-point and nine-point input vectors. Consider $\Delta x=\pi/4,\pi/8,\pi/16,\pi/32$ and $\Delta t=\Delta x=\Delta y$. In Table \ref{test1-table1}-\ref{test1-table2}, we list the $L_2$ and $L_{\infty}$ errors  at final time $T=\pi$.  Training error for five point scheme is shown in Fig. \ref{1-4errors} part (a).

\begin{table}[htbp]
\caption{Errors for five-point scheme with different initial values in Example \ref{ex1} }\label{test1-table1}
\begin{center}
\renewcommand{\arraystretch}{1.2}
\begin{tabular}{cccccc}
\hline
\multirow{2}{*}{$\Delta x$} &  \multicolumn{2}{c}{$cos(x+y)$} & &\multicolumn{2}{c}{$cos(x+y+\pi/3)$}\\
\cline{2-3}\cline{5-6}
&$L_2$ &$L_{\infty}$ &&$L_2$ &$L_{\infty}$\\
\hline
$\pi/4$ &   2.3250e-3 &4.8762e-4 & &2.5987e-3 &6.8945e-4\\
$\pi/8$ &9.2396e-4 &2.0430e-4 & &9.2396e-4 &2.2890e-4\\
$\pi/16$ &3.2872e-4 & 6.6806e-5 && 3.2872e-4 &6.6778e-5\\
$\pi/32$ & 1.9411e-4 &4.4076e-5 &&1.9276e-4 &4.4360e-5\\
\hline
\end{tabular}
\end{center}
\end{table}

\begin{figure}[htbp]
\centering
\subfigure{
\includegraphics[width=7cm]{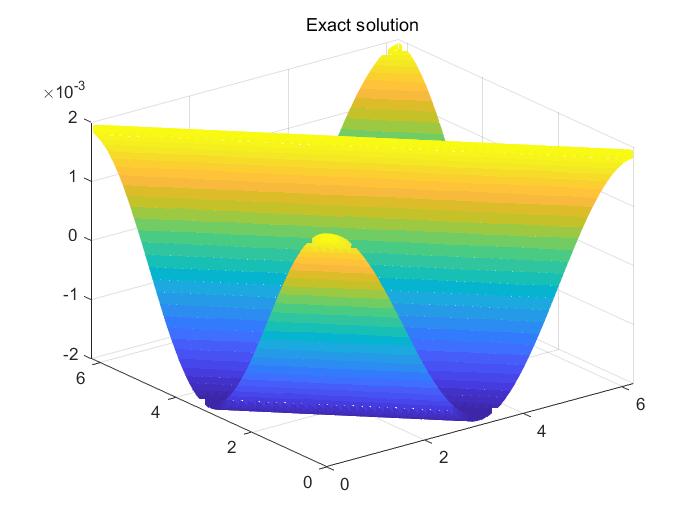} } \quad
\subfigure{ \includegraphics[width=7cm]{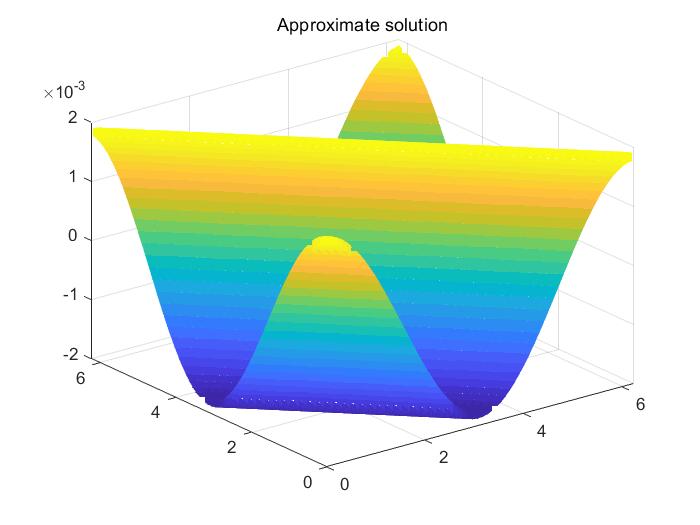} }
\caption{Exact solution(left) and approximate solution  with $\Delta x=\Delta t=\pi/32$ with initial value $\cos(x+y)$ }
\label{heat_figure}
\end{figure}

\begin{table}[htbp]
\caption{Errors for nine-point scheme with different initial values in Example \ref{ex1}}\label{test1-table2}
\begin{center}
\renewcommand{\arraystretch}{1.2}
\begin{tabular}{cccccc}
\hline
\multirow{2}{*}{$\Delta x$} &  \multicolumn{2}{c}{$cos(x+y)$} && \multicolumn{2}{c}{$cos(x+y+\pi/3)$}\\
\cline{2-3}\cline{5-6}
&$L_2$ &$L_{\infty}$ &&$L_2$ &$L_{\infty}$\\
\hline
$\pi/4$ & 2.3250e-3 &6.6374e-4 &&3.0197e-3 &7.2560e-4\\
$\pi/8$ &1.0971e-3 &2.3161e-4 &&1.0971e-4 &2.3534e-4\\
$\pi/16$ &3.2967e-4 & 6.5078e-5 && 3.2967e-4 &6.5082e-5\\
$\pi/32$ & 1.9154e-4 &4.3589e-5 &&1.9194e-4 &4.3174e-5\\
\hline
\end{tabular}
\end{center}
\end{table}

\vspace{.1in}
 \noindent {\bf Case II: test time step size.}

 We test the effect of time step size  on the accuracy and stability of CANN method. We choose five-point input vector and fix the spatial mesh size $\Delta x=\Delta y=\pi/16$, but vary the time sizes from $\Delta t=\Delta x/2, \Delta t=\Delta x, \Delta t=2\Delta x$ to $\Delta t=4\Delta x$. We compute the $L_2$ and $L_{\infty}$ errors and list in Table \ref{test1-table3}-\ref{test1-table4}. Numerical results indicate the accuracy and stability of CANN method seems to be independent of $\Delta t$.

 \begin{table}[htbp]
  \caption{Errors for five-point scheme with different initial value with $\Delta x=\pi/16$ for Example \ref{ex1}}\label{test1-table3}
 \begin{center}
\renewcommand{\arraystretch}{1.2}
\begin{tabular}{cccccc}
\hline
\multirow{2}{*}{$\Delta t$} &  \multicolumn{2}{c}{$cos(x+y)$} && \multicolumn{2}{c}{$cos(x+y+\pi/3)$}\\
\cline{2-3}\cline{5-6}
&$L_2$ &$L_{\infty}$ &&$L_2$ &$L_{\infty}$\\
\hline
$4\Delta x$ & 5.5740e-4 &1.2273e-4 &&5.6214e-4 &1.1715e-4\\
$2\Delta x$ &3.8369e-4 &7.8043e-5 &&3.8522e-4 &7.5607e-5\\
$\Delta x$ &3.2869e-4 & 6.8315e-5 && 2.9954e-4 &6.2361e-5\\
$\Delta x/2$ & 3.1106e-4 &6.4661e-5 &&2.9053e-4 &6.0928e-5\\
\hline
\end{tabular}
\end{center}
 \end{table}

\begin{table}[htbp]
\caption{Errors for nine-point scheme with different initial value with $\Delta x=\pi/16$ for Example \ref{ex1}}\label{test1-table4}
\begin{center}
\renewcommand{\arraystretch}{1.2}
\begin{tabular}{cccccc}
\hline
\multirow{2}{*}{$\Delta t$} &  \multicolumn{2}{c}{$cos(x+y)$} && \multicolumn{2}{c}{$cos(x+y+\pi/3)$}\\
\cline{2-3}\cline{5-6}
&$L_2$ &$L_{\infty}$ &&$L_2$ &$L_{\infty}$\\
\hline
$4\Delta x$ & 5.1559e-4 &1.1587e-4 &&5.2316e-4 &1.1160e-4\\
$2\Delta x$  &3.1420e-4 &6.1804e-5 &&3.2921e-4 &6.3846e-5\\
$\Delta x$ &3.1926e-4 & 6.0286e-5 && 3.1152e-4 &6.3144e-5\\
$\Delta x/2$& 2.9921e-4 &5.9587e-5 &&2.9830e-4 &6.0141e-5\\
\hline
\end{tabular}
\end{center}
\end{table}

%
%\vspace{.1in}
%\noindent {\bf Case III: test accuracy of five-point for different initial conditions.}
%
%For this case , we use the same solvers to solve the hear equation with the initial conditions $u=cos(x+y+\pi/3)$.
%In Table \ref{test1-table4}, we list the $L_2$ and $L_{\infty}$ errors and orders at final time $T=\pi$.
%\begin{table}[htbp]
%\caption{Errors and orders for CANN method with  five-point scheme in Example \ref{ex1} }\label{test1-table4}
%\centering
%\begin{tabular}[c]{l l l l l l}
%\hline
%$\Delta x$ & $L_2$ & order & $L_{\infty}$ & order\\
%\hline $\pi /4$  & 2.5987e-3 & & 6.8945e-4 &\\
%$\pi /8$  & 9.2396e-4 &1.49 & 2.2890e-4 &1.59\\
%$\pi/16$  & 3.2872e-4 &1.49 & 6.6778e-5  &1.77\\
%$\pi/32$  & 1.9276e-4 &0.77 & 4.4360e-5 &0.59\\
%\hline
%\end{tabular}
%
%\end{table}

  %%%%%%%%%%%%%%%%%%%%%%%%%%%%%%%%%%%%%%%%%%%%%
  \begin{example}\label{ex2}
  {\bf\emph {Convection diffusion equation}}
  \end{example}
  Consider the linear convection diffusion equation
  \begin{equation}\label{eq:3-2} u_t+c(u_x+u_y)-\mu(u_{xx}+u_{yy})=0,~(x,y)\in(0,2\pi)\times (0,2\pi), \end{equation}
  with the constants $c$, $\mu$  and initial condition $u(x,y,0)=\sin(x+y)$.  The training data is from the the problem with initial condition $u(x,y,0)=\sin(x+y)$ and the training network is applied to solve the same PDE with initial value  $u(x,y,0)=\cos(x+y)$ and $u(x,y,0)=\cos(x+y+\pi/6)$. The network picked 1 hidden layer and 15 neurons per layer. Iterations of $K=10^5$ is used to obtain the optimal parameter set. Total errors and orders are computed at the final time $T=\pi$.  Let $c=\mu=1$.

\vspace{.1in}
  \noindent {\bf Case I: test accuracy of nine-point scheme.}

  Four spatial mesh sizes of $\Delta x=\pi/4, \pi/8, \pi/16$ and $\pi/32$ are studied. We choose $\Delta t=\Delta x$ correspondingly. The well-trained network is used to solve the problem with both initial values $u(x,y,0)=\cos(x+y)$ and $u(x,y,0)=\cos(x+y+\pi/6)$. The corresponding errors are shown  in Table \ref{test2-table1} at final time $T=\pi$. Training error is shown in Fig. \ref{1-4errors} part (b).

   \begin{table}[htbp]
   \caption{Errors for different initial values in Example \ref{ex2}}\label{test2-table1}
  \begin{center}
\renewcommand{\arraystretch}{1.2}
\begin{tabular}{cccccc}
\hline
\multirow{2}{*}{$\Delta x$} &  \multicolumn{2}{c}{$cos(x+y)$} & & \multicolumn{2}{c}{$cos(x+y+\pi/6)$}\\
\cline{2-3}\cline{5-6}
&$L_2$ &$L_{\infty}$ & &$L_2$ &$L_{\infty}$\\
\hline
$\pi/4$ & 4.7801e-3 &1.0627e-3 &&4.8287e-3 &1.1093e-3\\
$\pi/8$ &1.6575e-3 &3.9201e-4 & &1.6576e-3 &3.9159e-4\\
$\pi/16$ &9.6424e-4 & 2.6655e-4 && 9.6524e-4 &2.6852e-4\\
$\pi/32$ & 3.1186e-4 &9.6218e-5 &&3.7686e-4 &1.1267e-4\\
\hline
\end{tabular}
\end{center}
  \end{table}

   %\begin{table}[htbp]
  % \caption{Errors and orders of CANN method in Example \ref{ex2}, %$T=\pi$}\label{test2-table4}
  %\centering
 % \begin{tabular}[c]{l l l l l l}
 % \hline $\Delta x$ & $L_2$ & order & $L_{\infty}$ & order\\
 % \hline $\pi /4$  & 4.8287e-3 & & 1.1093e-3 &\\
 % $\pi /8$  & 1.6575e-3 &1.54 & 3.9159e-4 &1.50\\
 % $\pi/16$  & 9.6524e-4 &0.78 & 2.6655e-4 &0.54\\
 % $\pi/32$  & 3.7686e-4 &1.36 & 1.1267e-4 &1.25\\
 % \hline \end{tabular}
 % \end{table}

\vspace{.1in}
  \noindent {\bf Case II: test time step size.}

  We choose nine-point input vector (\ref{nn:input-vector-stencil}) and fix the spatial mesh size $\Delta x=\pi/16$ with varying the time step size from $\Delta t=\Delta x/2, \Delta t=\Delta x, \Delta t=2\Delta x$ to $\Delta t=4\Delta x$. The errors for two initial values are shown in Table \ref{test2-table2}.
  \begin{table}[htbp]
  \caption{Errors for different initial values with $\Delta x=\pi/16$ for Example \ref{ex2}}\label{test2-table2}
  \begin{center}
\renewcommand{\arraystretch}{1.2}
\begin{tabular}{cccccc}
\hline
\multirow{2}{*}{$\Delta t$}  & \multicolumn{2}{c}{$cos(x+y)$} && \multicolumn{2}{c}{$cos(x+y+\pi/6)$}\\
\cline{2-3}\cline{5-6}
 &$L_2$ &$L_{\infty}$ &&$L_2$ &$L_{\infty}$\\
\hline
$4\Delta x$  & 1.5373e-4 &3.5473e-4& &1.5386e-3 &3.5533e-4\\
$2\Delta x$ &1.3389e-4 &3.5394e-4& &1.3389e-4 &3.5667e-4\\
$\Delta x$ &9.4351e-4 & 2.6616e-4 && 9.4351e-4 &2.6729e-4\\
$\Delta x/2$ & 6.6732e-4 &2.1398e-4 &&6.6732e-4 &2.1453e-4\\
\hline
\end{tabular}
\end{center}
  \end{table}

  %%%%%%%%%%%%%%%%%%%%%%%%%%%%%%%%%%%%%%%%%%%%%%%%%%%%%%%%%%%%%%%%%%%%%%%%%%%%%%%
  \begin{example}\label{ex3}
  {\bf\emph {Anisotropic diffusion equation}}
  \end{example}
  Consider anisotropic diffusion equation described as \begin{equation}\label{eq:3-3}
  u_t +c(u_x+u_y)= \mu(u_{xx}+u_{xy}+u_{yy}),  ~\text{in} ~  (0,2\pi)\times (0,2\pi)\times (0,T)
  \end{equation}
 with initial value $u(x,y,0)=\cos(x+y)$, periodic boundary condition and  final time $T=\pi/4$. Let $c=1$ and $\mu =0.01$. The training data is from the initial value $u(x,y,0)=\cos(x+y)$ and approximate solution at $t_1$ time level obtained from discontinuous Galerkin numerical method and then we use the well trained network to solve the same equation with different initial value $u(x,y,0)=\sin(x+y)$.  Take 1 hidden layer and 15 neurons per layer. The optimal parameter sets could be acquired after iterations $K=10^5$.

\vspace{.1in}
  \noindent {\bf Case I: test accuracy of five-point scheme.}

  We consider four setting of $\Delta x=\pi/16,\pi/32,\pi/64,\pi/128$ refined to check out the accuracy of the method while employ $\Delta t=\Delta x$. In Fig. \ref{1-4errors} part (c), we present the neural network squared $L_2$ training errors. Well-trained network errors are around $10^{-9}$, and the error hardly change  after the number of iterations $K=10^2$.  The errors and orders over $L_2$ and $L_{\infty}$ for different initial values at $T=\pi/4$ are shown in Table \ref{test3-table1}.

  \begin{table}[!htb]
     \caption{Errors for initial values $u(x,y,0)=sin(x+y)$ in Example \ref{ex3}}\label{test3-table1}
\begin{center}
\begin{tabular}{ccccc}
\hline
$\Delta x$ & $L_2$&order &$L_{\infty}$ &order \\
\hline
$\pi/16$ & 7.3895e-4 & &1.8229e-4&  \\
$\pi/32$ &6.1777e-4&0.26 & 1.5343e-4 &0.25 \\
$\pi/64$ &8.6468e-5&2.84 & 1.9292e-5 &2.99 \\
$\pi/128$ &3.8696e-5& 1.16 &6.7268e-6 &1.52 \\
\hline
\end{tabular}
\end{center}
\end{table}

\vspace{.1in}
   \noindent {\bf Case II: test time step size.}

   In this test, we change the time size.  While we fix the spatial size $\Delta x=\pi/32$ and adjust $\Delta t$ for four times, as $\Delta t=4\Delta x,~ 2\Delta x,~\Delta x,~ \Delta x/2$. The  errors for different initial values are listed in Table \ref{test3-table2} at final time $T=\pi$.
   \begin{table}[htbp]
   \caption{Errors for initial values $u(x,y,0)=sin(x+y)$ with $\Delta x=\pi/32$ for Example \ref{ex3}}\label{test3-table2}
   \begin{center}

\begin{tabular}{ccc}
\hline
&$L_2$ &$L_{\infty}$ \\
\hline
$4\Delta x$ & 7.9598e-4 &2.1803e-4 \\
$2\Delta x$ &4.8973e-4 & 1.2276e-5 \\
$\Delta x$ &3.2068e-4 & 6.9930e-5  \\
$\Delta x/2$ &4.1742e-4 &1.1344e-4  \\
\hline
\end{tabular}
\end{center}
   \end{table}
   %%%%%%%%%%%%%%%%%%%%%%%%%%%%%%% %%%%%%%%%%%%%%%%%%%%%%%%%%%%%%

   \begin{example}
   \label{ex4}
   {\bf\emph {2D incompressible Navier-Stokes equation in vorticity formulation}} \end{example}
   Consider 2D incompressible Navier-Stokes equations in vorticity  formulation
   \begin{equation}
   \label{eq:3-4}
   u_t+\nabla \cdot (wu)=\frac{1}{Re} \Delta u,~ (x,y)\in (0,2\pi)\times (0,2\pi),
   \end{equation}
   with periodic boundary condition and initial value $u(x,y,0)=2\cos y\sin x$. To simplify the computation, choose incompressible velocity field $w=(-\cos x\sin y, \sin x\cos y)$. The work of \cite{jin2021nsfnets} motivate us to study this with employing the Reynolds number $Re=100$. In the following tests, we choose the input vector (\ref{nn:input-vector-stencil}), 1 hidden layer, 15 neurons and the number of iteration $K=10^5$. The training data is obtain from the problem with the initial condition $u(x,y,0)=2\cos y\sin x$ and test the same PDE with initial condition  $u(x,y,0)=2\cos x\sin y$.

\vspace{.1in}
   \noindent {\bf Case I: test accuracy of five-point scheme.}

\begin{figure}[htbp]
   \centering
   \subfigure{ \includegraphics[width=7cm]{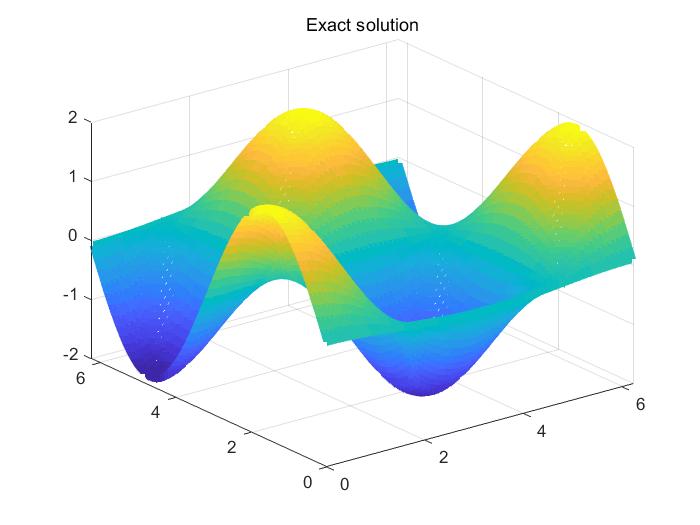}} \quad \subfigure{ \includegraphics[width=7cm]{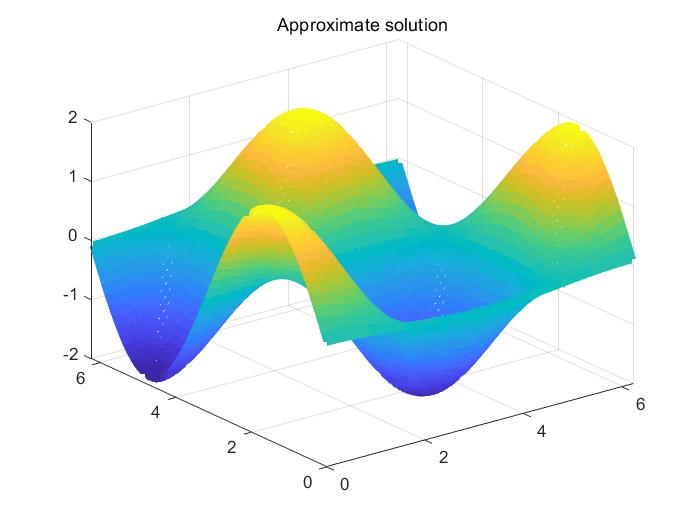} }
   \caption{Exact solution(left) and approximate solution(right) at  $T=\pi$ in Example \ref{ex4} with $\Delta x=\Delta t=\pi/32$.}\label{Navier_figure}
   \end{figure}

  Convergence experiments are undertaken first by using four setting of $\Delta x=\pi/4, \pi/8,\pi/16,\pi/32$. For each cell size, we apply $\Delta x=\Delta y=\Delta t$. In Fig. \ref{1-4errors} part (d), we conclude that we could acquire the optimal parameter after iterations of almost $K=10^5$. $L_2$ and $L_{\infty}$ errors of the CANN method are calculated at final time $T=\pi$ and listed in Table \ref{test4-table1}. Again we show the comparison of solution of the incompressible Navier-Stokes equation in Fig.\ref{Navier_figure}.

   \begin{table}[htbp]
   \caption{Errors for different initial value $u(x,y,0)=2\cos x\sin y$ in Example \ref{ex4}}\label{test4-table1}
  \begin{center}
\begin{tabular}{ccc}
\hline

&$L_2$ &$L_{\infty}$  \\
\hline
$\pi/4$ & 4.3961e-4 &9.0003e-5  \\
$\pi/8$ &3.7362e-4 & 1.4306e-4  \\
$\pi/16$ &9.9108e-5 & 4.3985e-5 \\
$\pi/32$ &1.5989e-4 &7.3565e-5  \\
\hline
\end{tabular}
\end{center}
   \end{table}

\vspace{.1in}
   \noindent {\bf Case II: test time step size.}

   In this case, the spatial mesh size is fixed as $\Delta x=\pi/16$ but we apply four time sizes $\Delta t=4\Delta x, 2\Delta x, \Delta x$ and $\Delta x/2$. The corresponding results are in Table \ref{test4-table2}. The accuracy show that CANN method seems to be independent of time step size $\Delta t$.
   \begin{table}[!htb]
   \caption{Errors for different initial value with $\Delta x=\pi/16$ for Example \ref{ex4}}\label{test4-table2}
  \begin{center}

\begin{tabular}{ccc}
\hline

&$L_2$ &$L_{\infty}$ \\
\hline
$4\Delta x$ & 6.3322e-4 &2.8864e-4  \\
$2\Delta x$ &4.2085e-4 & 1.8836e-4 \\
$\Delta x$ &1.2048e-4 & 2.8652e-4 \\
$\Delta x/2$ &1.0750e-4 &7.8368e-5 \\
\hline
\end{tabular}
\end{center}
   \end{table}

 %%%%%%%%%%%%%%%%%%%%%%%%%%%%%%%%%%%%%%%%%%%%%%%%%%%%%%%%%%%%
    \begin{example}
    \label{ex5}{\bf\emph {Nonlinear  equation}}
    \end{example}
    In this section, we discuss the following nonlinear problem: \begin{equation}\label{eq:3-5}
    u_t  =\nabla \cdot a(x,\nabla u) + f,   ~ \text{in} ~  \Omega \times (0,T),\quad \text{with}\ \ u(x,y,0)=(x^2+y^2)/2,
    \end{equation}
     where $\Omega=[-1,1]\times [-1,1]$ and $T=1$. Specifically, we take $a(x,\nabla u)=(1+e^{-\| \nabla u \|^2})\nabla u$, then we have $f=(4e^{-2t-r}-1)u-2e^{-t}(e^{-r}+1)$ which $r=-2e^{-t}u$. The exact solution is $u=e^{-t}(x^2+y^2)/2$. In the tests, we take five-point scheme as the input vector,  1 hidden layer, 15 neurons per layer, the number of iterations $K=10^5$ and final time $T=1$.
We randomly select $75\%$ of solution average data (\ref{nn:learning-data-set}) with $p=0$ at  $t_0$ and $t_1$ time level as the training set and the other $25\%$ is used as the testing data. %Then we obtain the approximate solution at $t_n, n\geq 2$ time level from the trained network.

\vspace{.1in}
     \noindent {\bf Case I: test accuracy of five-point scheme.}

     Four spatial meshes of $\Delta x=1/4, 1/8, 1/16, 1/32$ are considered and let $\Delta x=\Delta t$ to explore the convergence of solvers which use CANN method to approximate nonlinear diffusion equation. The errors of the testing data at $t_1$ and  $T=1$ time level are shown in Table \ref{test5-table1}. % rrors and orders at $T=1$ in Table \ref{test7-table1}.

     \begin{table}[htbp]
     \caption{ Errors  for the testing data at $t_1$ and $T=1$ in Example \ref{ex5}}\label{test5-table1}
     \begin{center}
\renewcommand{\arraystretch}{1.2}
\begin{tabular}{cccccc}
\hline
\multirow{2}{*}{$\Delta x$} & \multicolumn{2}{c}{testing data at $t_1$} && \multicolumn{2}{c}{testing data at $T=1$} \\
\cline{2-3}\cline{5-6}
&$L_2$ &$L_{\infty}$& &$L_2$ &$L_{\infty}$ \\
\hline
$1/4$ & 1.8262e-3 &3.0337e-3 && 1.2983e-2 &1.0371e-2 \\
$1/8$ &8.7103e-4 & 2.0438e-3 &&7.1876e-3 &1.0109e-2 \\
$1/16$ &1.9588e-4 & 1.0610e-3 &&3.2560e-3 & 7.4244e-3 \\
$1/32$ &5.0547e-5 &3.6399e-4 && 1.2614e-3 &4.4760e-3 \\
\hline
\end{tabular}
\end{center}
     \end{table}

   %  \begin{table}[htbp]
    % \caption{Errors and orders of CANN methods for Example \ref{ex7}.}\label{test7-table1}
    % \centering
    % \begin{tabular}[c]{l l l l l l}
    % \hline $\Delta x$ & $L_2$ & order & $L_{\infty}$ & order\\
    % \hline $1/4$  & 1.4268e-2 & & 1.1413e-2&\\
    % $1/8$  & 7.5624e-3 &0.92 & 1.0506e-2 &0.12\\
    % $1/16$  & 3.3076e-3 &1.19 & 7.6508e-3&0.46\\
    % $1/32$  & 1.3649e-3 &1.28 & 4.6466e-3 &0.72\\
    % \hline
    % \end{tabular}
    % \end{table}

      %\begin{figure}[htbp]
     % \centering \includegraphics[width=8.5cm]{figure/3-1-5.jpg} \caption{Squared $L_2$ training errors (\ref{nn:L2-error-square}) for Example \ref{eq:3-7}.}\label{nonlinear_figure}
     % \end{figure}

      \vspace{.1in}
     \noindent {\bf Case II: test time step size.}

     In this case, we demonstrate that our method can keep validity in nonlinear problem. The spatial cell size $\Delta x=1/16$ is fixed. Again considering the four time meshes of $\Delta t=4\Delta x, 2\Delta x, \Delta x$ and $\Delta x/2$. We list all $L_2$ and $L_{\infty}$ errors in Table \ref{test5-table2}. The four well-trained CANN solvers  give the similar errors. Although the equation becomes complex, the method seems allow large time size for neural network solver.
     \begin{table}[htbp]
     \caption{Errors  for the testing data at $t_1$ and $T=1$ with $\Delta x=1/16$ in Example \ref{ex5}}
     \label{test5-table2}
     \begin{center}
\renewcommand{\arraystretch}{1.2}
\begin{tabular}{cccccc}
\hline
\multirow{2}{*}{$\Delta t$} & \multicolumn{2}{c}{testing data at $t_1$} && \multicolumn{2}{c}{testing data at $T=1$} \\
\cline{2-3}\cline{5-6}
&$L_2$ &$L_{\infty}$ &&$L_2$ &$L_{\infty}$ \\
\hline
$4\Delta x$ & 6.7075e-4 &3.5711e-3 && 2.2281e-3 &5.2523e-3 \\
$2\Delta x$ &4.1082e-4 & 2.1802e-3&&2.8491e-3 &6.4610e-3 \\
$\Delta x$ &2.3186e-4 & 1.1720e-3 &&3.3055e-3 & 7.4138e-3 \\
$\Delta x/2$ &1.0605e-5 &5.6397e-4 && 3.5455e-3 &8.0069e-3 \\
\hline
\end{tabular}
\end{center}
     \end{table}

     %%%%%%%%%%%%%%%%%%%%%%%%%%%%%%%%%%%%%%%
     \begin{example}\label{ex6}
     {\bf\emph {Porous medium equation}}
     \end{example}
     In this section, we consider the two-dimensional porous medium equation (PME) under the bounded spatial region $\Omega=[0,1]\times [0,1]$ and time interval $t=[0,1]$ which can be defined as \cite{chew2017newton}:
     \begin{equation}\label{eq:3-6}
     u_t =0.2 \nabla\cdot(u^2\nabla u),  ~ \text{in} ~  \Omega \times (0,T),\quad \text{with}\ \ u(x,y,0)=(5(x+y)+15)^{1/2}.
          \end{equation}
     The solution is $u(x,y,t)=(5(x+y+t)+15)^{1/2}$. In the tests, we take five-point scheme as the input vector,  1 hidden layer, 6 neurons per layer, the number of iterations $K=10^5$ and final time $T=1$. We generated the training data from the problem (\ref{eq:3-6}) and then test
     the same PDE with different initial value $u(x,y,0)=(5(x+y)+11)^{1/2}$.

\vspace{.1in}
     \noindent {\bf Case I: test accuracy of five-point scheme.}

\begin{figure}[htbp]
   \centering \subfigure[Example \ref{ex5}: Nonlinear equation]{
   \includegraphics[width=6.5cm]{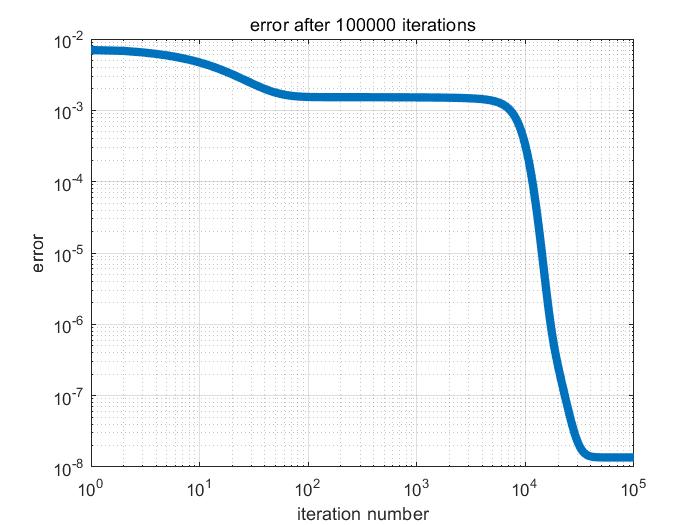} } \quad
   \subfigure[Example \ref{ex6}: Porous medium equation]{
   \includegraphics[width=6.5cm]{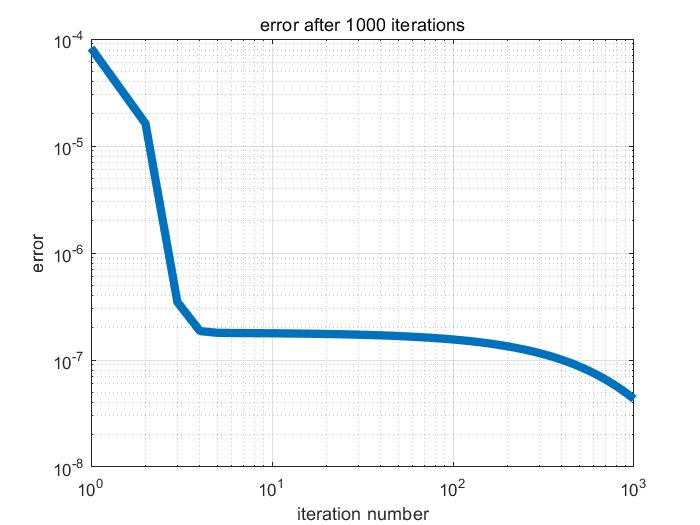} } \quad
   \subfigure[Example \ref{ex7}: 3D diffusion equation]{
   \includegraphics[width=6.5cm]{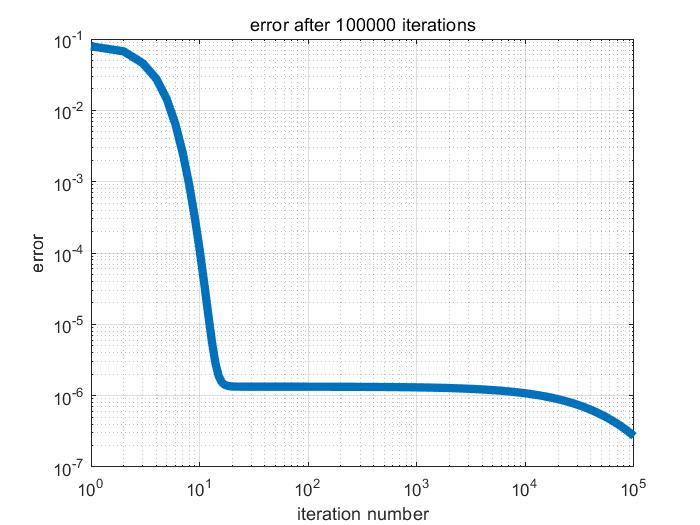} } \quad
   \subfigure[Example \ref{ex8}: 4D diffusion equation]{
   \includegraphics[width=6.5cm]{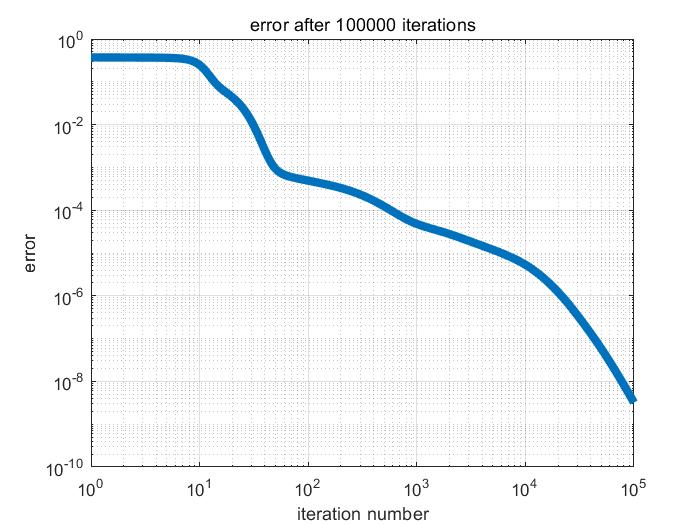} } \caption{Neural network squared $L_2$ training errors
   (\ref{nn:L2-error-square})}\label{2-4errors}
   \end{figure}

     Consider four setting of $\Delta x=\Delta y=1/4, 1/8, 1/16,1/32$ refined to check out the convergence and let $\Delta t=\Delta x$. In Fig. \ref{2-4errors} part(b) , we show the neural network squared $L_2$ training errors. Well-trained network errors are around $10^{-8}$. In Table \ref{test6-table1}, $L_2$ and $L_{\infty}$ errors at $T=1$ for initial value $u(x,y,0)=(5(x+y)+11)^{1/2}$ is listed.
    % \begin{figure}[htbp]
     %\centering
    % \includegraphics[width=8.5cm]{figure/PME_error.jpg}
    % \caption{Squared $L_2$ training errors (\ref{nn:L2-error-square}) for Example %\ref{eq:3-8}.}\label{PME_figure}
    % \end{figure}
     \begin{table}[htbp]
     \caption{Errors for different initial value $u(x,y,0)=(5(x+y)+11)^{1/2}$ in Example \ref{ex6}}\label{test6-table1}
    \begin{center}

\begin{tabular}{ccccc}
\hline
&$L_2$&order&$L_{\infty}$ &order\\
\hline
$1/4$ & 8.4116e-3& &2.3618e-2& \\
$1/8$ &3.5175e-3 & 1.26& 1.4228e-2&0.74  \\
$1/16$ &3.2383e-3&0.12 & 1.0758e-2 & 0.40\\
$1/32$ &3.2335e-3& 0.01 &9.9734e-3&0.11  \\
\hline
\end{tabular}
\end{center}
     \end{table}

    \vspace{.1in}
     \noindent {\bf Case II: test time step size.}

     In this part, the spatial mesh size is fixed as $\Delta x=1/16$ but we take four time size $\Delta t=4\Delta x,2\Delta x,\Delta x,\Delta x/2$. The computed $L_2$ and $L_{\infty}$ errors of the four network solvers are listed in Table \ref{test6-table2}. Notice the choice of $\Delta t=4\Delta x$ involves a time size roughly 128 times bigger than the regular CFL restriction of $\Delta t=O(\Delta x^2)$.
     \begin{table}[htbp]
     \caption{Errors for different initial value $u(x,y,0)=(5(x+y)+11)^{1/2}$ with $\Delta x=1/16$ for Example \ref{ex6}}\label{test6-table2}
     \begin{center}

\begin{tabular}{ccc}
\hline

&$L_2$ &$L_{\infty}$ \\
\hline
$4\Delta x$ & 3.5421e-3 &1.6366e-2  \\
$2\Delta x$ &3.5341e-3 & 1.7164e-2 \\
$\Delta x$ &3.1066e-3 & 9.5569e-2  \\
$\Delta x/2$ &3.5729e-3 &1.0799e-2 \\
\hline
\end{tabular}
\end{center}
     \end{table}

     %%%%%%%%%%%%%%%%%%%%%%%%%%%%%%%%%%%%%%%%%%%%%%%%
     \subsection{Higher dimensional diffusion equations} \label{S:3.2}

     In this section, we investigate the higher dimensional diffusion equations. From two-dimensional problem in Section \ref{S:3.1}, we have demonstrate that the approximate solutions for our CANN method have high accuracy for large time step size and the convergence of the method seems be independent of the CFL restrictions. Suppose that the computed domain is a hypercube and every cell is still a hypercube after subdivision.  We take the cell average of the solution in each cell $\Gamma$ and the cell average in all the cells adjacent to $\Gamma$ as the input vectors of the neural network similar to the scheme (\ref{nn:co-vertex-input-stencil}).   Given the dimension $d$, every cell has $3^d$ cells adjacent to it. We use two order finite difference method to obtain the approximate solution at $t=t_1=\Delta t$ as the output vector of the neural network. We take $\Delta x_1=\Delta x_2=\cdots =\Delta x_d$ and the number of iteration $K=10^5$. In the following test, we consider the  problem
     \begin{equation}\label{eq:4}
     \begin{aligned}
     u_t&= \Delta u, & ~ \text{in} ~  \Omega \times (0,T),\\
     u &=u|_{t=0}, & \text{in} ~ \Omega \times \{0\},\\
     u &=u_0, ~& \text{on} ~\partial \Omega \times (0,T),
     \end{aligned}
     \end{equation}
     with the domain $\Omega=[0,\pi]^d$, the dimension of the equation $d\geq 3$ and periodic boundary condition $u_0$.  Take $T=\pi$ as the final time.

     %%%%%%%%%%%%%%%%%%%%%%%%%%%%%%%%%%%%%%%%%%%%%%%%%%%%%%%%%%%%%
     \begin{example}\label{ex7}
     {\bf\emph {3D diffusion equation}}
     \end{example}

     Consider 3D diffusion equation (\ref{eq:4}) with periodic boundary and initial condition $u(x_1,x_2,x_3,0)=\sin(x_1+x_2+x_3)$. We can give the exact solution is $u(x_1,x_2,x_3,t)=e^{-3t}cos(x_1+x_2+x_3)$. Especially, we get training data by the equation (\ref{eq:4}) with initial condition $u(x_1,x_2,x_3,0)=\sin(x_1+x_2+x_3)$ and test the same equation with different initial value $u(x_1,x_2,x_3,0)=\cos(x_1+x_2+x_3)$.
     We use the similar network structure in Section \ref{S:2.1}-\ref{S:2.4}. The input vector can be defined as \begin{equation}\label{eq:3-D-co-vertex-input-vector}
     \vv_{ijk}^{in} = \Big[ \uAve^n_{i-1,j-1,k-1},\uAve^n_{i-1,j-1,k},\uAve^n_{i-1,j-1,k+1},\cdots,\uAve^n_{i-1,j+1,k+1},\uAve^n_{i,j+1,k+1},\uAve^n_{i+1,j+1,k+1}\Big]^T.
     \end{equation}
      The dimension of input vector is $3^d=27$. We take 1 hidden layer and 10 neurons per layer. The relation between $L_2$ training errors and the number of iteration is shown in Fig. \ref{2-4errors} part(c). Test the accuracy of CANN method with the four meshes $\Delta x_i=\pi/2^{i+1}, i=1,2,3,4$ and the effect of time step size on the convergence of the method. At final time $T=\pi$, The corresponding results are shown in Table \ref{test7-table1}.
      \begin{table}[htbp]
      \caption{Errors  for different initial value $\Delta t=\Delta x_i$ in Example \ref{ex7}}\label{test7-table1}
      \begin{center}

\begin{tabular}{ccccc}
\hline

&$L_2$&order &$L_{\infty}$&order \\
\hline
$\pi/4$ & 3.5384e-3& &6.6643e-3& \\
$\pi/8$ &2.6042e-3 & 0.44&4.8890e-3&0.45 \\
$\pi/16$ &3.2450e-4&3.00 & 6.0500e-4&3.01 \\
$\pi/32$ &1.0321e-4&1.65 &1.9683e-4&1.62 \\
\hline
\end{tabular}
\end{center}
      \end{table}

      %\begin{figure}[htbp]
      %\centering \includegraphics[width=8.5cm]{figure/3d_error.jpg} \caption{training errors and the number of iteration for  Example \ref{ex9} with $\Delta x_i=\pi/32$.}\label{3d_figure}
      %\end{figure}

      \begin{table}[htbp]
      \caption{Errors for different initial value with $\Delta x=\pi/16$ for Example \ref{ex7}}\label{test7-table2}
     \begin{center}

\begin{tabular}{ccc}
\hline

&$L_2$ &$L_{\infty}$ \\
\hline
$4\Delta x$ &  4.7848e-4 &5.3145e-4 \\
$2\Delta x$ &3.7432e-4 &6.8006e-4 \\
$\Delta x$ &3.4556e-4 & 6.1667e-4 \\

\hline
\end{tabular}
\end{center}
      \end{table}

      %%%%%%%%%%%%%%%%%%%%%%%%%%%%%%%%%%%%%%%%%%%%%%%%%%%%%%%%%%%%%%%%%%%%%%%%%%%%%%

      \begin{example}\label{ex8}
      {\bf\emph {4D diffusion equation}}
      \end{example}

      In this section, we consider the problem (\ref{eq:4}) with $d=4$, periodic boundary condition and initial condition $u(x_1,x_2,x_3,x_4,0)=\sin(x_1+x_2+x_3+x_4)$.    Take 1 hidden layer and 10 neurons every layer and final time $T=\pi$. The input vector of the network is similar to (\ref{eq:3-D-co-vertex-input-vector}). The training data of (\ref{nn:learning-data-set}) is generated by the equation (\ref{eq:4}) with initial condition $u=(x_1,x_2,x_3,x_4,t)=exp(-4t)\sin(x_1+x_2+x_3+x_4)$ and the trained network is applied to test the same PDE with different initial value $u(x_1,x_2,x_3,x_4,0)=\cos(x_1+x_2+x_3+x_4)$.
      The dimension of the input vector is $81$.
      Similarly, we test the accuracy of the method and the effect of time step size on the accuracy and the corresponding results are shown in Fig. \ref{2-4errors} part(d) and in Table \ref{test8-table1}.
     % \begin{figure}[htbp]
      %\centering
      %\includegraphics[width=8.5cm]{figure/4d_error.jpg} \caption{training errors and the number of iteration for  Example \ref{ex10} with $\Delta x_i=\pi/16$. }\label{4d_figure}
     % \end{figure}

      \begin{table}[!htb]
      \caption{Errors for different initial value  $\Delta t=\Delta x_i$ in Example \ref{ex8}}\label{test8-table1}
     \begin{center}

\begin{tabular}{ccccc}
\hline

&$L_2$&order &$L_{\infty}$&order \\
\hline
$\pi/4$ & 4.8038e-3& &8.7656e-4 &\\
$\pi/8$ &2.2889e-4 & 4.39 &2.5155e-5&5.12 \\
$\pi/16$ &3.1760e-5& 2.85 & 4.7712e-6 &2.40\\

\hline
\end{tabular}
\end{center}
      \end{table}

 \section{Conclusions}  \label{S:4}

      In this article, we develop a CANN method to solve high-dimensional  parabolic problems. First, we use numerical methods to obtain the approximate solution at $t=t_1=\Delta t$ as the output of the network. Secondly, we take the initial value as input vectors, get the well-trained network and test the same PDE with different initial value. Finally, we compute the approximate solution at $t=t_{n+1}$  for given solution at $t=t_n$ with the trained optimal network and present the error for the approximate solution at $t=T$ time level. All the numerical results in Sections \ref{S:3.1}-\ref{S:3.2} indicate the following conclusions:

      1) CANN method is convergent and has remarkable accuracy for large time step size.

      2) CANN method is found to be relieved from the strict CFL condition also in high dimensional problem and could evolve the solution forward in time with large time step size.
     % When we use numerical methods to compute the approximate solution of the original problem at time $t=\Delta t=t_1$, spatial step size and time step size must met the CFL condition.
%      Later we use the initial value as input vector and the solution at $t=t_1$ as output vector to train the neural network and obtain optimal network parameters. We use the trained network to  compute the approximate solution at $t_{n+1}$ from the solution at $t_n, n\geq 1$ and at this time  time step size $\Delta t$ may not be restricted to CFL conditions.

      3) The trained network can be employed to solve the problem with the same equation with the different initial condition.

       %%%% Acknowledgments %%%%%%%%
       %\section*{Acknowledgments}

     %  The research of Huang was supported by National Natural Science Foundation of China under grant NO.11771398.  Research work of Yan is supported by National Science Foundation grant DMS-1620335 and Simons Foundation grant 637716.

       %\bibliographystyle{elsarticle-num}

       %\bibliography{cas-refs}

\begin{thebibliography}{99}
\bibitem{qiu2021cell}
C.~Qiu and J.~Yan, Cell-average based neural network method for hyperbolic and parabolic partial differential equations, arXiv preprint arXiv:2107.00813.

\bibitem{bengio2009learning}
Y.~Bengio, Learning deep architectures for AI, Now Publishers Inc., 2009.

\bibitem{beck2019machine}
C.~Beck, E.~Weinan and A.~Jentzen, Machine learning approximation algorithms for
  high-dimensional fully nonlinear partial differential equations and
  second-order backward stochastic differential equations, Journal of Nonlinear
  Science, 29~(4) (2019), 1563-1619.

\bibitem{goodfellow2016deep}
I.~Goodfellow, Y.~Bengio and A.~Courville, Deep learning, MIT press, 2016.

\bibitem{lecun2015deep}
Y.~LeCun, Y.~Bengio and G.~Hinton, Deep learning, Nature, 521~(7553) (2015),
  436-444.

\bibitem{rackauckas2020universal}
C.~Rackauckas, Y.~Ma, J.~Martensen, C.~Warner, K.~Zubov, R.~Supekar,
  D.~Skinner, A.~Ramadhan and A.~Edelman, Universal differential equations for
  scientific machine learning, arXiv preprint arXiv:2001.04385.

\bibitem{chen2018neural}
R.~T. Chen, Y.~Rubanova, J.~Bettencourt and D.~Duvenaud, Neural ordinary
  differential equations, arXiv preprint arXiv:1806.07366.

\bibitem{long2019pde}
Z.~Long, Y.~Lu and B.~Dong, PDE-Net 2.0: Learning PDEs from data with a
  numeric-symbolic hybrid deep network, Journal of Computational Physics, 399
  (2019), 108925.

\bibitem{ruthotto2020deep}
L.~Ruthotto and E.~Haber, Deep neural networks motivated by partial differential
  equations, Journal of Mathematical Imaging and Vision, 62~(3) (2020), 352-364.

  \bibitem{he2019mgnet}
J.~He and J.~Xu, Mgnet: A unified framework of multigrid and convolutional neural
  network, Science china mathematics, 62~(7) (2019) 1331-1354.

\bibitem{ray2018artificial}
D.~Ray and J.~S. Hesthaven, An artificial neural network as a troubled-cell
  indicator, Journal of computational physics, 367 (2018), 166-191.

\bibitem{wang2019learning}
Y.~Wang, Z.~Shen, Z.~Long and B.~Dong, Learning to discretize: solving 1d scalar
  conservation laws via deep reinforcement learning, arXiv preprint
  arXiv:1905.11079.

\bibitem{discacciati2020controlling}
N.~Discacciati, J.~S. Hesthaven and D.~Ray, Controlling oscillations in high-order
  discontinuous galerkin schemes using artificial viscosity tuned by neural
  networks, Journal of Computational Physics, 409 (2020), 109304.

\bibitem{lagaris1998artificial}
I.~E. Lagaris, A.~Likas and D.~I. Fotiadis, Artificial neural networks for solving
  ordinary and partial differential equations, IEEE transactions on neural
  networks, 9~(5) (1998), 987-1000.

\bibitem{rudd2015constrained}
K.~Rudd and S.~Ferrari, A constrained integration (cint) approach to solving
  partial differential equations using artificial neural networks,
  Neurocomputing, 155 (2015), 277-285.

\bibitem{sirignano2018dgm}
J.~Sirignano and K.~Spiliopoulos, DGM: A deep learning algorithm for solving
  partial differential equations, Journal of computational physics, 375 (2018),
  1339-1364.

\bibitem{raissi2017physics}
M.~Raissi, P.~Perdikaris and G.~E. Karniadakis, Physics informed deep learning
  (part I): Data-driven solutions of nonlinear partial differential equations,
  arXiv preprint arXiv:1711.10561.

\bibitem{raissi2019physics}
M.~Raissi, P.~Perdikaris and G.~E. Karniadakis, Physics-informed neural networks:
  A deep learning framework for solving forward and inverse problems involving
  nonlinear partial differential equations, Journal of Computational Physics,
  378 (2019), 686-707.

\bibitem{dwivedi2020physics}
V.~Dwivedi and B.~Srinivasan, Physics informed extreme learning machine (pielm)--a
  rapid method for the numerical solution of partial differential equations,
  Neurocomputing, 391 (2020), 96-118.

\bibitem{lu2021deepxde}
L.~Lu, X.~Meng, Z.~Mao and G.~E. Karniadakis, DeepXDE: A deep learning library for
  solving differential equations, SIAM Review, 63~(1) (2021), 208-228.

\bibitem{jin2021nsfnets}
X.~Jin, S.~Cai, H.~Li and G.~E. Karniadakis, Nsfnets (Navier-Stokes flow nets):
  Physics-informed neural networks for the incompressible navier-stokes
  equations, Journal of Computational Physics, 426 (2021), 109951.

\bibitem{shin2020convergence}
Y.~Shin, J.~Darbon and G.~E. Karniadakis, On the convergence of physics informed
  neural networks for linear second-order elliptic and parabolic type PDEs,
  arXiv preprint arXiv:2004.01806.

\bibitem{laakmann2021efficient}
F.~Laakmann and P.~Petersen, Efficient approximation of solutions of parametric
  linear transport equations by ReLU DNNs, Advances in Computational
  Mathematics, 47~(1) (2021), 1-32.

\bibitem{cai2021least}
Z.~Cai, J.~Chen and M.~Liu, Least-squares relu neural network (LSNN) method for
  linear advection-reaction equation, Journal of Computational Physics, (2021),
  110514.

\bibitem{cai2021least1}
Z.~Cai, J.~Chen and M.~Liu, Least-squares relu neural network (LSNN) method for
  scalar nonlinear hyperbolic conservation law, arXiv preprint
  arXiv:2105.11627.

\bibitem{Weinan2018Ritz}
W.~E and B.~Yu, The deep Ritz method: a deep learning-based numerical algorithm
  for solving variational problems, Communications in mathematics and
  statistics, 6 (2018), 1-12.

\bibitem{zang2020weak}
Y.~Zang, G.~Bao, X.~Ye and H.~Zhou, Weak adversarial networks for high-dimensional
  partial differential equations, Journal of Computational Physics, 411 (2020),
  109409.

\bibitem{chan2019machine}
Q.~Chan-Wai-Nam, J.~Mikael and X.~Warin, Machine learning for semilinear PDEs,
  Journal of Scientific Computing, 79~(3) (2019), 1667-1712.

\bibitem{pham2021neural}
H.~Pham, X.~Warin and M.~Germain, Neural networks-based backward scheme for fully
  nonlinear PDEs, SN Partial Differential Equations and Applications, 2~(1)
  (2021), 1-24.

\bibitem{hutzenthaler2020proof}
M.~Hutzenthaler, A.~Jentzen, T.~Kruse and T.~A. Nguyen, A proof that rectified
  deep neural networks overcome the curse of dimensionality in the numerical
  approximation of semilinear heat equations, SN partial differential equations
  and applications, 1~(2) (2020), 1-34.

\bibitem{lye2020deep}
K.~O. Lye, S.~Mishra and D.~Ray, Deep learning observables in computational fluid
  dynamics, Journal of Computational Physics, 410 (2020), 109339.

\bibitem{wu2020data}
K.~Wu and D.~Xiu, Data-driven deep learning of partial differential equations in
  modal space, Journal of Computational Physics, 408 (2020), 109307.

\bibitem{qin2021data}
T.~Qin, Z.~Chen, J.~D. Jakeman and D.~Xiu, Data-driven learning of nonautonomous
  systems, SIAM Journal on Scientific Computing, 43~(3) (2021), A1607--A1624.

\bibitem{chew2017newton}
J.~Chew and J.~Sulaiman, Newton-sor iterative method for solving the
  two-dimensional porous medium equation, Journal of Fundamental and Applied
  Sciences, 9~(6S) (2017), 384-394.

     \end{thebibliography}

       \end{document}